\documentclass[12pt]{amsart}
\usepackage{amsmath,amsthm,latexsym,amscd,amsbsy,amssymb,url}
\usepackage[all]{xypic}
 \relax

\chardef\bslash=`\\ 

\makeatletter
\def\verbatim{\interlinepenalty\@M \@verbatim
  \leftskip\@totalleftmargin\advance\leftskip2pc
  \frenchspacing\@vobeyspaces \@xverbatim}
\makeatother
\newtheorem{thm}{Theorem}[section]

\newtheorem{lem}[thm]{Lemma}

\newtheorem{defin}[thm]{Definition}
\newtheorem{rem}[thm]{Remark}

\begin{document}


\title
{On Topological Groups of Automorphisms on Unions}
\author{Raushan  Buzyakova}
\email{Raushan\_Buzyakova@yahoo.com}



\begin{abstract}{
We study groups of homeomorphic bijections on spaces that are finite unions of compact connected  linearly ordered subsets. We prove that all such groups when endowed with the topology of point-wise convergence are topological groups.
}
\end{abstract}

\maketitle
\markboth{R. Buzyakova}{On Topological Groups of Monotonic Automorphisms}
{ }

\section{Introduction}\label{S:introduction}
\par\bigskip\noindent
In this paper, a homeomorphic bijection of a topological space with itself will be called an  automorphism.
We will study  group operations on   the space   $Hom_p(X)$ of automorphisms on $X$ endowed with the topology of point-wise convergence. When we refer to $Hom_p(X)$ as a topological group we mean that the operation of composition is continuous on $Hom_p(X)\times Hom_p(X)$ and the map $f\mapsto f^{-1}$ is continuous on $Hom_p(X)$. We will be  concerned with the following general problem:

\par\bigskip\noindent
{\bf Problem} {\it
Let $Hom_p(X)$ and $Hom_p(Y)$  be topological groups.

\begin{enumerate}
	\item Is $Hom_p(X\times \{0,1\})$ a topological group?
	\item Is $Hom_p(X\oplus Y)$ a topological group?
	\item Is  $Hom_p(X\cup Y)$ a topological group, where $X$ and $Y$ are closed subspaces of the union?
\end{enumerate}
}

\par\bigskip\noindent
It was proved in \cite{Sor} that $Hom_p(L)$ is a topological group for any connected LOTS $L$.  In this work we give affirmative answers to all three questions of the problem if $X$ and $Y$ are connected and compact LOTS by proving a more general statement, namely that $Hom_p(X)$ is a topological group whenever $X$ is the union of a finite collection of connected compact LOTS.

In notation and terminology of the general topological nature, we will follow \cite{Eng}. A space is always a topological space. Recall that  an open subset $U$ of a space $X$ is canonical if $U$ is the interior of $\overline U$. For general facts about topological groups, we refer the reader to \cite{AT}. A linearly ordered topological spaces (or simply ordered space), abbreviated as LOTS, is a linearly ordered set endowed with the topology generated by sets $(a,b)$, $\{x\in L: x<a\}$, and $\{x\in L: x>a\}$. For general LOTS-related facts, we refer the reader to \cite{L}. Finally, when dealing with several linearly ordered sets at the same time, we  will distinguish their intervals via subscription as in $[a,b]_L$ (the same applies to other types of intervals). 

\par\bigskip
\section{Preliminary Observations}

\par\bigskip\noindent
In this section, we will discuss the topological structure of spaces that  can be written as the union of a finite number of compact connected linearly ordered topological subspaces. Given such a space $X$, we will identify a basis for the topology of $X$ that will be used in the argument of our main result.

\begin{defin}
A topological space $X$ is locally ordered at $x\in X$ if there exists an open neighborhood of $x$ which is a LOTS. The set of all points at which $X$ is locally ordered is denoted by $LO(X)$.
\end{defin}

\par\bigskip\noindent
\begin{lem}\label{lem:LO}
Let $X$ be the union of a finite number of closed linearly ordered spaces. Then $LO(X)$ is an open and dense subspace of $X$.
\end{lem}
\begin{proof}
Let $X = \bigcup\{L_i: i =1,...,n\}$, where  $L_i$ is a closed subset of $X$ and is a LOTS for each $i\in I$.
Fix any non-empty set $U\subset X$. We need to find a non-empty open subset of $U$ which is a LOTS. Let $k\leq n$ be the smallest such that $U$ is a subset of 
 $\bigcup\{L_i: i =1,...,k\}$.  Therefore, $V= U\setminus  \bigcup\{L_i: i =1,...,k-1\}$ is not empty and is open in $X$. Since $V\subset L_k$, V is an a non-empty open subset of $U$ which is a LOTS.
\end{proof}

\par\bigskip\noindent
\begin{defin}\label{defin:orderly}
Let $X$ be the union of a finite number of compact connected LOTS. An open set $U\subset X$ is called orderly if the following properties hold:
\begin{enumerate}
	\item $U$ is canonical.
	\item $U$ is connected.
	\item $\overline U\setminus U$ is finite.
	\item$\overline U\setminus U\subset LO(X)$
\end{enumerate}
\end{defin}

\par\bigskip\noindent
\begin{lem}\label{lem:twoorderlysets}
Let $U$ and $V$ be connected open subsets of $X$ and $\overline U\setminus U= \overline V\setminus V$. Then either $U=V$ or $U\cap V=\emptyset$.
\end{lem}
\begin{proof}
Let $x\in V\cap U$. It suffices to show that $V\subset U$. Assume that $V\setminus U$ is not empty. Since $V$ is connected, $V\setminus U$ is not open. Then there exists $y\in V\setminus U$ such that any neighborhood of $y$ meets both $U$ and $V$. Since $y\not \in U$, we conclude that $y\in \overline U\setminus U$. Since $\overline U\setminus U= \overline V\setminus V$, we arrive at $y\not \in V$, a contradiction.
\end{proof}

\par\bigskip\noindent
The following Lemma \ref{lem:randomopendense} will be referenced  both for its argument and for the statement.
\par\bigskip\noindent
\begin{lem}\label{lem:randomopendense}
Let $X$ be the union of a finite number  of compact connected LOTS.
Let $U$ be an orderly neighborhood of $p\in X$, and let $D$ be a dense subset of $X$. Then there exists an orderly neighborhood $V$ of $p$ such that 
$\overline V \subset U$ and $\overline V \setminus V\subset D$.
\end{lem}
\begin{proof} If $p\in LO(X)$, then the conclusion is clear. We now assume that $p\not\in LO(X)$. Put $B=\overline U\setminus U=\{b_i: i =1,...,n\}$. For each $k\in \{1,...,n\}$, let $I_k\subset LO(X)$ be  an open neighborhood of $b_k$ such that $\overline I_k \cap \overline I_{m}=\emptyset$ for distinct $k,m$. 
Since $U$ is canonical, we may assume that $I_k$ is a connected LOTS without both maximum and minimum. We may also assume that $(b_k, \max I_k]_{I_k}\subset U$. For each $b_k\in B$, 
fix $c_k\in (b_k, \max I_k)_{I_k}\cap D$.  Let $V_0 = U$ and $V_k = V_{k-1}\setminus (b_k, c_{k}]_{I_k}$ for $0<k\leq n$. Put $V= V_n$. Let us show that $V$ is as desired.
Since $p\not \in LO(X)$, we conclude that $p\in V$. Since $(b_k,c_k]_{I_k}$ is closed in $U$ for each $b_k\in B$, the set $V_k$ is open for each $k$. 

Next let us show that $V$ is orderly. The boundary of $V_k$ is   $\{b_i: i=1,...,k-1\}\cup \{c_i:i=k,...,n\}$.This follows from the fact that $c_k$ is an interior point of $I_k$ for each $k$. Therefore, $\overline V\setminus V=\{c_k: k=1,...,n\}$ is a subset of $D\cap LO(X)$. It remains to show that $V$ is connected and canonical.

To show that $V$ is canonical, observe that  any neighborhood of  $c_k$ has a smaller neighborhood that is a subset of $I_k$ and is of the form $(x,y)_{I_k}$. Therefore, $(c_k, y)_{I_k}$ and $(x,c_k)_{I_k}$ are non-empty sets that meet $V$ and the complement of $V$, respectively.

To show that $V$ is connected, let us show that $V_k$ is connected for each $k=0,...,n$. Since $U$ is orderly, $V_0$ is connected.
 Assume that $V_k$ is connected for $0\leq k <n$. To show that $V_{k+1}$ is connected let us argue by contradiction. Then, 
there exist disjoint open sets $O$ and $O'$ such that $V_{k+1} = O \cup O'$ . Then $V_k\setminus \{c_{k+1}\}$ is the union of disjoint open sets  $(b_{k+1}, c_{k+1})_{I_{k+1}}$, $O$, and $O'$. Since $V_{k+1}$ is connected $c_{k+1}$ is a boundary point of each of the three sets contradicting the fact that $c_{k+1}$ has an open neighborhood which is a connected LOTS. 
\par\medskip\noindent
Our proof is complete and let us finish with a remark for future use.
\par\bigskip\noindent
{\it Remark. Using  the constructions of $V$ we can create an orderly set $W$  by taking $c_k$ anywhere in $I_k$ not necessarily to the right of $b_k$ and by letting
$W= U \cup\bigcup\{I_k\setminus (min I_k, c_k]_{I_k}: k=1,...,n\}$. This set is also orderly by a similar argument but  need not be a subset of $U$ if at least one  $c_k$ is taken outside of $U$. 
}
\end{proof}

\par\bigskip\noindent
\begin{rem}\label{rem:construction} The sets described in the final remark of  Lemma \ref{lem:randomopendense} will be useful in our future arguments. Let us refer to the described $W$ as being \underline{determined by $U, \{I_k\}_k, \{c_k\}$}.
\end{rem}

\par\bigskip\noindent
If in the statement of In Lemma \ref{lem:randomopendense}, we put additional restrictions on $D$ and the boundary of $U$, the same argument leads to the following statement.

\par\bigskip\noindent
\begin{lem}\label{lem:randomopendense1}
Let $X$ be the union of a finite number of compact connected LOTS.
Let $U$ be an orderly neighborhood of $p\in X$, and let $D$ be an open dense subset of $X$.  Suppose that $\overline U\setminus U\subset D$. 
Then there exists an orderly neighborhood $V$ of $p$ such that 
$\overline V \subset U$ and $U \setminus V\subset D$.
\end{lem}

\par\bigskip\noindent
\begin{lem}\label{lem:orderlybasis}
Let $X$ be the union of a finite number of compact connected LOTS. Then $X$ has a basis consisting of orderly open sets.
\end{lem}
\begin{proof}
Let $X = \bigcup_{i=1}^nL_i$, where $L_i$ is a connected compact LOTS. We will induct on $n$. If $n=1$, then $X$ is a LOTS and the conclusion follows. Assume now that the conclusion of our lemma holds for $n=k$, where $k\geq 1$. We now assume that $n=k+1$. Fix $p\in X$ and an open neighborhood $O$ of  $p$.  Our goal is  to find an orderly neighborhood of $p$ which is a subset of $O$. Let $X' = \bigcup_{i=1}^{n-1}L_i$ and $O'=O\cap X'$. We may assume that $p\in X'$. By our assumption there exists an orderly open neighborhood $U$ of $p$ in $X'$ which is a subset of $O'$. 

\par\bigskip\noindent
{\it Claim 1. $LO(X)\cap O'$ is open and dense in $O'$.}
\par\smallskip\noindent
To prove the claim, it suffices to show that the interior  $Int_X(X')$ of $X'$ in $X$ is dense in $X'$. Pick $x\in X'\setminus Int_X(X')$. Since $L_n$ is closed in $X$, we have $x\in L_n$. Since  $x\not\in Int_X(X')$, we have $x\not \in LO(X)$. Therefore, any neighborhood of $x$ contais points that are not in $L_n$. Since $L_n$ is closed , all such points are in the interior of $X'$ in $X$.
\par\bigskip\noindent
By Claim 1 and Lemmas 
\ref{lem:LO} and \ref{lem:randomopendense}, we may assume that $\overline U\setminus U\subset LO(X)$. By Lemma \ref{lem:randomopendense1}, we can find an orderly neighborhood $V$ of $p$ in $X'$ such that $\overline V\subset U$ and $U\setminus V\subset LO(X)$. 
\par\smallskip\noindent
Now let us color the points of $G = \overline V$ in green color and the points of $R= X\setminus U$ in red color.
\par\bigskip\noindent
{\it Claim 2. There exists a  finite collection $\mathcal J$  of convex open sets in $L_n$ such that the following hold:
\begin{enumerate}
	\item $\bigcup\mathcal J$ is a subset of $O$.
	\item $\bigcup\mathcal J$ contains all green points of $L_n$, and misses all red points of $L_n$.
	\item Each $I\in \mathcal J$ has green points.
	\item The closures of distinct members of $\mathcal J$ are disjoint.
\end{enumerate}
}
\par\smallskip\noindent
The conclusion of the claim follows from the fact that the sets of green  and red points in $L_n$ are disjoint compacta and all green points are in $O$.

\par\bigskip\noindent
{\it Remark: Due to density and openness of $LO(X)$ in $X$, we may assume that for each 
$J\in \mathcal J$, the set $\overline J\setminus J \subset LO(X)$.}
\par\medskip\noindent
Put $W = U\cup \bigcup{\mathcal J}$. Let us show that $W$ is an orderly neighborhood of $p$ and is a subset of $O$. Since $U$ contains $p$ so does $W$. Since both $U$ and $\bigcup {\mathcal J}$ are subsets of $O$ so is $W$. To prove that $W$ is open, fix an arbitrary $x\in W$. We have two cases:
\begin{description}
	\item[$x\in L_n$]  Let $U_x$ be an open neighborhood of $x$ in $X$ that misses $R$ (all read points of $X$) as well as $L_n\setminus \bigcup\mathcal J$. Then $U_x\cap L_n\subset \bigcup\mathcal J$ and  $U_x\cap X'\subset U$. Hence, $U_x\subset W$.
	\item[$x\not \in L_n$] Then there exists a neighborhood $U_x$ of $x$ in $X$ that misses $L_n$ and $X'\setminus U$. Hence $U_x\subset W$.
\end{description}
The set $W$ is connected since $U$ and  $J$'s are and each $J$ contains some green points of $U$. The boundary of $W$ is finite since  every point on the boundary of $W$ is either on the boundary of $U$ or on the boundary of some $J\in \mathcal J$, which are finite. Finally, if $W$ is not canonical we can simply replace it with the interior of $\overline W$.
The proof is complete.
\end{proof}

\par\bigskip
\section{Study}

\par\bigskip\noindent
We are now ready to prove our main result that if $X$ is the union of a finite number of connected LOTS, then $Hom_p(X)$ is a topological group. For this we need to verify that the operation of function composition is a continuous map from $Hom_p(X)\times Hom_p(X)$ to $Hom_p(X)$ and that the correspondence $f\mapsto f^{-1}$ is a continuous map from $Hom_p(X)$ to $Hom_p(X)$. 

\par\bigskip\noindent
\begin{lem}\label{lem:inversion}
Let $X$ be the union of a finite number of connected compact LOTS. Then $f\mapsto f^{-1}$ is a continuous map from $Hom_p(X)$
to $Hom_p(X)$.
\end{lem}
\begin{proof}
Fix $f\in Hom_p(X)$. Let $V_{f^{-1}}$ be any open neighborhood of $f^{-1}$. We need to find an open $U_f$ containing $f$ such that $g^{-1}\in V_{f^{-1}}$ whenever $g\in U_f$.
We may assume that there exist $x\in X$ and an orderly neighborhood $O_x$ of $x$ such that $V_{f^{-1}} = \{h\in Hom_p(X): h(y)\in O_x\}$, where $y=f(x)$. Let $O'_x$ be an orderly neighborhood of $x$ such that $\overline {O_x'}\subset O_x$. For each $b\in \overline {O'_x}\setminus O'_x$, fix an open neighborhood $I_b$ of $b$ in $LO(X)\cap O_x$ such that $x\not \in \overline I_b$ and  $\{\overline I_b: b\in \overline {O'_x}\setminus O'_x\}$ is disjoint. Since $O_x$ is canonical, we may assume that $I_b$ is a connected LOTS without both maximum and minimum. Let $O''_x$ be an open neighborhood of $x$  that misses each $I_b$. Let
$$
U_f =\{h\in Hom_p(X): h(b)\in f(I_b)\ for\ every\ b\in \overline {O'_x}\setminus O'_x\ and\ h(x) \in f(O_x'') \}.
$$
Since the border of $O'_x$ is finite and the images of open sets under $f$ are open, the set is an open neighborhood of $f$. To show that $U_f$ is as desired, fix an arbitrary $h\in U_f$. 

\par\bigskip\noindent
{\it Claim. 
$h(O'_x)$ contains $y$.
}
\par\smallskip\noindent
To prove the claim first note that $h(O'_x)$ is an orderly set containing $h(x)$. According to Remark \ref{rem:construction}, the orderly set $B$ determined by $f(O'_x), \{f(I_b) : b\in \overline {O'_x}\setminus O'_x\}$, and $\{h(b) : b\in \overline {O'_x}\setminus O'_x\}$ also contains $h(x)$. By Lemma \ref{lem:twoorderlysets}, the sets coincides. By Remark  \ref{rem:construction},  $B$ contains $y$ since $f(I_b)$ misses $f(x)=y$ for each $b$ in the border of $O'_x$.

\par\bigskip\noindent
By Claim, $h^{-1}(y)\in O'_x\subset O_x$. Hence, $h^{-1}\in V_{f^{-1}}$.
\end{proof}

\par\bigskip\noindent
\begin{lem}\label{lem:key}
Let $X$ be the union of a finite number of compact connected LOTS. Let $f\in Hom_p(X)$ and $U_y$ an open neighborhood of $y=f(x)$. Then there exist open neighborhoods $U_x$ and $U_f$ of $x$ and $f$, respectively such that $h(U_x)\subset U_y$ whenever $h\in U_f$
\end{lem}
\begin{proof}
Let $O_x$ be an orderly neighborhood of $x$ such that $f(O_x)\subset U_y$.  Let $O_x'$ be an orderly neighborhood of  $x$ such that $\overline {O'_x}\subset O_x$. For each $b\in \overline {O'_x}\setminus O'_x$, fix an open neighborhood $I_b$ of $b$ in $LO(X)\cap O_x$  such that $x\not \in \overline I_b$ and $\{\overline I_b: b\in \overline {O'_x}\setminus O'_x\}$ is disjoint. We can assume that each $I_b$ is a connected LOTS without  both maximum and minimum. Let
$$
U_f =\{h\in Hom_p(X): h(b)\in f(I_b)\ for\ every\ b\in \overline {O'_x}\setminus O'_x\ and\ h(x) \in f(O_x') \}.
$$
\par\medskip\noindent
Since our construction is very similar to that in Lemma \ref{lem:inversion}, 
by the argument identical to that of  Claim of Lemma \ref{lem:inversion}, $h(O'_x)$ contains $f(a)$ for every $a\in O_x'$ which is not in 
$\cup\{\overline I_b: b\in \overline {O'_x}\setminus O'_x\}$ whenever $h\in U_f$. Therefore, $U_x=O'_x\setminus \cup\{\overline I_b: b\in \overline {O'_x}\setminus O'_x\}$ is as desired.
\end{proof}

\par\bigskip\noindent
\begin{lem}\label{lem:composition}
Let $X$ be the union of a finite number of compact connected LOTS. Then $\langle f, g\rangle \mapsto f\circ g$ is a continuous map  from $Hom_p(X)\times Hom_p(X)$ to $Hom_p(X)$.
\end{lem}
\begin{proof}
Pick $f, g\in Hom_p(X)$ and an open neighborhood $W_{f\circ g}$ of  $f\circ g$. We need to find open neighborhoods $U_f$ and $V_g$ of $f$ and $g$, respectively, such that $f'\circ g'\in W_{f\circ g}$ whenever $f'\in U_f$ and $g'\in V_g$.We may assume that there exists an open set $O_z\subset X$ of $z=f((g(x))$ such that  $W_{f\circ g} = \{h\in Hom_p(X): h(x)\in O_z\}$.

By Lemma \ref{lem:key}, there exist a neighborhood $U_f$ of $f$  and a neighborhood $O_y$ of $y$ such that $h(O_y)\subset O_z$ whenever $h\in U_f$.
Next, By Lemma \ref{lem:key}, there exist a neighborhood $V_g$ of $g$  and a neighborhood $O_x$ of $x$ such that $h(O_x)\subset O_y$ whenever $h\in V_f$. Clearly, $U_f$ and $V_g$ are as desired.
\end{proof}

\par\bigskip\noindent
Lemmas \ref{lem:inversion} and  \ref{lem:composition} imply our main result.

\par\bigskip\noindent
\begin{thm}
Let $X$ be the union of a finite number of connected compact LOTS. Then, $Hom_p(X)$ is a topological group.
\end{thm}
\par

\end{document}